\documentclass[journal]{IEEEtran}
\usepackage[cmex10]{amsmath}
\usepackage{latexsym,amsmath,enumerate,amssymb,amsbsy}
\usepackage{verbatim}
\usepackage{color}
\usepackage{url}
\usepackage{graphicx}
\newtheorem{theorem}{Theorem}[section]
\newtheorem{lemma}[theorem]{Lemma}

\newtheorem{definition}[theorem]{Definition}

\newtheorem{proposition}[theorem]{Proposition}

\newtheorem{example}[theorem]{Example}

\newcommand{\ZZ}{\mathbb{Z}}

\title{A Classification of Unimodular Lattice Wiretap Codes in Small Dimensions}
\author{Fuchun~Lin~and~Fr\'ed\'erique~Oggier
\thanks{The authors are with the Division of Mathematical Sciences,
  School of Physical and Mathematical Sciences, Nanyang Technological
  University, 21 Nanyang Link, Singapore 637371 
  (emails:linf0007@e.ntu.edu.sg and frederique@ntu.edu.sg). Part of this work appeared at ITW 2011\cite{ITWversion}.}%

\thanks{The research of F. Lin and of F. Oggier for this work is supported by the Singapore National Research
Foundation under the Research Grant NRF-RF2009-07.}%
}

\begin{document}

\maketitle

\begin{abstract}
Lattice coding over a Gaussian wiretap channel, where an eavesdropper listens to transmissions between a transmitter and a legitimate receiver, is considered. A new lattice invariant called the secrecy gain \cite{ISITA}  is used as a code design criterion for wiretap lattice codes since it was shown to characterize the confusion that a chosen lattice can cause at the eavesdropper: the higher the secrecy gain of the lattice, the more confusion. In this paper, a formula for the secrecy gain of unimodular lattices is derived. Secrecy gains of extremal odd unimodular lattices as well as unimodular lattices in dimension $n,\ 16\leq n \leq 23$ are computed, covering the $4$ extremal odd unimodular lattices and all the $111$ non-extremal unimodular lattices (both odd and even) providing thus a classification of the best wiretap lattice codes coming from unimodular lattices in dimension $n,\ 8<n \leq 23$. Finally, to permit lattice encoding via Construction A, the corresponding error correction codes are determined.
\end{abstract}



\section{Introduction}\label{sec:introduction}
In his seminal work, Wyner [2] introduced the wiretap channel, a discrete memoryless channel where the sender Alice
transmits confidential messages to a legitimate receiver Bob, in the presence of an eavesdropper Eve.
Both reliable and confidential communication between Alice and Bob is shown to be achievable at the same time,
by exploiting the physical difference between the channel to Bob and that to Eve, without the use of cryptographic means. Many results of information theoretical nature are available in the literature for various classes of channels ranging from Gaussian point-to-point channels to relay networks (see e.g. \cite{theoretic security} for a survey) capturing the trade-off between reliability and secrecy and aiming at determining the highest information rate that can be achieved with perfect secrecy, the so-called \textit{secrecy capacity}. Coding results focusing on constructing concrete codes that can be implemented in a specific channel are much fewer (see \cite{OW-84,TDCMM-07} for examples of wiretap codes dealing with channels with erasures).

In this paper, we will focus on Gaussian wiretap channels, whose secrecy capacity was established in \cite{IEEE Gaussian wiretap channel}. Examples of existing Gaussian wiretap codes were designed for binary inputs, as in \cite{LDPC Gaussian wiretap codes,nested codes}. A different approach was adopted in \cite{ISITA}, where lattice codes were proposed, using as design
criterion a new lattice invariant called \textit{secrecy gain}, which was shown to characterize the confusion at the eavesdropper. This suggests the study of the secrecy gain of lattices as a way to understand how to design a good Gaussian wiretap code. \textit{Unimodular} lattices were shown to be good candidates in \cite{ITW} and for \textit{even unimodular} lattices, both secrecy gains for a special class of lattices called \textit{extremal} lattices were computed and the asymptotic behavior of the average secrecy gain as a function of the dimension $n$ was investigated. These two papers were further developed in \cite{Gaussian wiretap codes}, where coding examples were detailed and it was shown that as $n$ grows to infinity, all even unimodular lattices behave in the same way, so that optimizing the secrecy gain makes sense in small dimensions.

The work of \cite{ITW,Gaussian wiretap codes} deals with even unimodular lattices, which only exist in dimensions a multiple of 8. We pursue the study of unimodular lattices by considering odd unimodular lattices, which on the contrary exist in every dimension and in great number, giving thus more flexibility in the code design. We will also show examples of odd unimodular lattices outperforming even unimodular lattices. Our contributions can be summarized as follows:
\begin{itemize}
\item We develop a general formula for the secrecy gain of both odd and even unimodular lattices that generalizes the existing one for even unimodular lattices.
\item We obtain the secrecy gain of unimodular lattices in dimension $n,\ 8<n\leq 23$, covering the $4$ extremal odd unimodular lattices as well as all the $111$ non-extremal unimodular lattices.
\item We classify the best Gaussian wiretap codes from unimodular lattices in dimension $n,\ 8<n\leq 23$, together with their corresponding self-dual codes enabling lattice encoding via Construction A.
\end{itemize}

The remainder of this paper is organized as follows. In Section II, we first give a brief introduction to unimodular lattices and their \textit{theta series} as well as recall the definition of the secrecy gain and the previous results concerning this lattice invariant. The main results are given in Section III. An explicit formula for the secrecy gain of unimodular lattices is derived, which generalizes the one for the even case in \cite{Gaussian wiretap codes}. Secrecy gains of extremal odd unimodular lattices are computed to complete the study of extremal unimodular lattices. Finally, secrecy gains of unimodular lattices in dimension $16\leq n\leq 23$, both odd and even, are computed, ending the classification of unimodular wiretap lattice codes in dimension $n,\ 8<n\leq 23$. In Section IV, encoding of the best codes via Construction A is discussed.


\section{Preliminaries and previous results}\label{sec:preliminaries}
Consider a Gaussian wiretap channel, which is modeled as follows: Alice wants to send data to Bob on a Gaussian channel whose noise variance is given by $\sigma_b^2$. Eve is the eavesdropper trying to intercept data through another Gaussian channel with noise variance $\sigma_e^2$, where $\sigma_b^2< \sigma_e^2$, in order to have a positive secrecy capacity \cite{IEEE Gaussian wiretap channel}. More precisely, the model is
\begin{equation}\label{channel model}
\begin{array}{cc}
\mathbf{y}&=\mathbf{x}+\mathbf{v_b}\\
\mathbf{z}&=\mathbf{x}+\mathbf{v_e},\\
\end{array}
\end{equation}
where $\mathbf{x}$ is the transmitted signal, $\mathbf{v_b}$ and $\mathbf{v_e}$ denote the Gaussian noise vectors at Bob's, respectively Eve's side, each component of both vectors with zero mean, and respective variance $\sigma_b^2$ and $\sigma_e^2$, and finally $\mathbf{y}$ and $\mathbf{z}$ are the received signals at Bob's, respectively Eve's side. In this paper, we choose $\mathbf{x}$ to be a codeword coming from a specially designed lattice of dimension $n$, namely, we consider lattice coding. Let us thus start by recalling some concepts concerning lattices, in particular, unimodular lattices.
\subsection{Unimodular lattices}

A \textit{lattice} $\Lambda$ is a discrete set of points in $\mathbb{R}^n$, which can be described in terms of its \textit{generator matrix} $M$ by 
$$
\Lambda=\{\mathbf{x}=\mathbf{u}M|\mathbf{u}\in \mathbb{Z}^n\}, 
$$  
where
$$
M=\left (
\begin{array}{cccc}
v_{11}&v_{12}&\cdots&v_{1n}\\
v_{21}&v_{22}&\cdots&v_{2n}\\
\cdots&&\cdots&\\
v_{n1}&v_{n2}&\cdots&v_{nn}\\
\end{array}
\right )
$$
and the row vectors $\mathbf{v}_i=(v_{i1},\cdots,v_{in}),\ i=1,\ 2,\ \cdots,\ n$ form a basis of the lattice. The matrix
$$
A=MM^T,
$$
where $M^T$ denotes the transpose of $M$, is called the \textit{Gram matrix} of the lattice. It is easy to see that the ($i$,$j$)th entry of $A$ is the inner product of the $i$th and $j$th row vectors of $M$, denoted by
$$A_{(i,j)}=\mathbf{v}_i\cdot \mathbf{v}_j.$$ 
A lattice $\Lambda$ is called an \textit{integral lattice} if its Gram matrix is an integral matrix. The \textit{determinant} det$\Lambda$ of a lattice $\Lambda$ is the determinant of the matrix $A$, which is independent of the choice of the matrix $M$. A \textit{fundamental region} for a lattice is a building block which when repeated many times fills the whole space with just one lattice point in each copy. There are many different ways of choosing a fundamental region for a lattice $\Lambda$, but the volume of the fundamental region is uniquely determined by $\Lambda$ and called the volume of $\Lambda$, which is exactly $\sqrt{\mbox{det}\Lambda}$. Let us see an example of a fundamental region of a lattice. A \textit{Voronoi cell} $\mathcal{V}_{\Lambda}(\mathbf{x})$ of a lattice point $\mathbf{x}$ in $\Lambda$ consists of the points in the space that are closer to $\mathbf{x}$ than to any other lattice points of $\Lambda$.

The \textit{dual} of a lattice $\Lambda$ of dimension $n$ is defined to be
$$
\Lambda^*=\{\mathbf{x}\in \mathbb{R}^n: \mathbf{x}\cdot \mathbf{\lambda}\in \mathbb{Z} , \mathbf{\lambda}\in \Lambda\}.
$$
It can be shown that $\Lambda$ is an integral lattice if and only if $\Lambda\subset \Lambda^*$. Especially, if $\Lambda=\Lambda^*$ then $\Lambda$ is called a \textit{unimodular} lattice. It can further be shown that $\Lambda$ is a unimodular lattice if and only if $\Lambda$ is integral and $\mbox{det }\Lambda=1$. Finally, the \textit{norm} (squared length) $||\mathbf{x}||^2=\mathbf{x}\cdot \mathbf{x}$ of a lattice point $\mathbf{x}$ in a unimodular lattice $\Lambda$ is of course an integer. If the norm is an even integer for any lattice point in $\Lambda$, then $\Lambda$ is called an \textit{even unimodular} lattice or a \textit{type II} lattice. Otherwise, it is called an \textit{odd unimodular} lattice or a \textit{type I} lattice. 

There are certain lattices which play the role of building blocks in analyzing lattices. They are denoted by $A_n$, $D_n$, $E_6$, $E_7$, $E_8$ and are called \textit{irreducible root lattices} \cite{lattices and codes}. Unimodular lattices are then decomposed into a number of such lattices. More precisely, an $n$-dimensional unimodular lattice $\Lambda$ is described as one containing a sublattice which is the direct sum
$$
\Lambda_1\oplus \Lambda_2\oplus \cdots \oplus \Lambda_k
$$
of a number of irreducible root lattices of total dimension $n$, and consequently a lattice point of $\Lambda$ can be written as
$$
\mathbf{x}=\mathbf{x_1}+\mathbf{x_2}+\cdots+\mathbf{x_k},
$$
where each component $\mathbf{x_i}$ is chosen as one of a standard system of representatives for the cosets of $\Lambda_i$ in $\Lambda_i^{*}$ and called a \textit{glue vector} for $\Lambda_i$. Informally speaking, $\Lambda$ is obtained by gluing together the components $\Lambda_1,\ \Lambda_2,\ \cdots,\ \Lambda_k$ by the glue vectors. The existence of glue vectors will be indicated by ``$^+$'' all through this paper. Let us see an example of the simplest case that is called \textit{self-glue}.

\begin{example}\label{ex:D_12+}
We use $(r^{m})$ to denote a string of $m$ $r$'s here. The lattice $D_{12}^+=D_{12}\cup D_{12}+(\frac{1}{2}^{12})$ contains $D_{12}$ as a sublattice and a lattice point of $D_{12}^+$ can be written as $\mathbf{x}=\mathbf{x_1},\ \mathbf{x_1}\in D_{12}+(0^{12})$ or  $\mathbf{x}=\mathbf{x_1},\ \mathbf{x_1}\in D_{12}+(\frac{1}{2}^{12})$, where the vectors $(0^{12})$ and $(\frac{1}{2}^{12})$ are the glue vectors for $D_{12}$.
\end{example}

A complete list of unimodular lattices of dimension $n,\ 0\leq n\leq 23$ that contain no vector of norm $1$ is given in \cite{sphere packings}, each lattice is described by its components. The \textit{kissing number} of a lattice (sphere) packing is the number of spheres that touch one sphere. The kissing numbers of these unimodular lattices are also given in the same table.

Let us recall the definition of the \textit{Jacobi theta functions} and the \textit{theta series} of lattices before we end this subsection. Let $\mathcal{H}=\{a+ib\in \mathbb{C}|b>0\}$ denote the upper half plane and let $q=e^{\pi i \tau}$, where $\tau\in \mathcal{H}$.
\begin{definition}
Jacobi theta functions are defined as follows: 
$$
\begin{array}{ll}
\vartheta_2(\tau)&=\Sigma_{n\in\mathbb{Z}}q^{(n+\frac{1}{2})^2},\\
\vartheta_3(\tau)&=\Sigma_{n\in\mathbb{Z}}q^{n^2},\\
\vartheta_4(\tau)&=\Sigma_{n\in\mathbb{Z}}(-q)^{n^2}.\\
\end{array}
$$
\end{definition}

They are very important functions in analytic number theory. For example, the \textit{discriminant function $\Delta_8(\tau)$} can be represented by $\vartheta_2(\tau)$ and $\vartheta_4(\tau)$ \cite{sphere packings}:
\begin{equation}\label{equ:discriminant}
\begin{array}{ll}
\Delta_8(\tau)&=\frac{1}{16}\vartheta_2^4(\tau)\vartheta_4^4(\tau)\\
                      &=2^{-8}\vartheta_2^8(\frac{\tau+1}{2})\\
                      &=q\prod_{m=1}^{\infty}\{(1-q^{2m-1})(1-q^{4m})\}^8\\
                      &=q-8q^2+28q^3-64q^4+126q^5+224q^6+\cdots.\\
\end{array}
\end{equation}

\begin{definition}
The theta series of a lattice $\Lambda$ is defined by
$$\Theta_{\Lambda}(\tau)=\Sigma_{\mathbf{\lambda}\in \Lambda}q^{\mathbf{\lambda}\cdot \mathbf{\lambda}}.$$
\end{definition}

If we combine the terms with the same exponent, the theta series of an integral lattice $\Lambda$ can be written as
\begin{equation}\label{equ:theta series}
\Theta_{\Lambda}(\tau)=\Sigma_{n=0}^{\infty}A_nq^n.
\end{equation}
By doing that we can interpret the theta series of $\Lambda$ as a book keeping device recording the number of vectors $\lambda\in \Lambda$ with norm $n$ in the coefficient $A_n$. Take the one-dimensional lattice $\mathbb{Z}$ for example and recall the definition of the Jacobi function $\vartheta_3(\tau)$. We have
\begin{equation}\label{equ:theta3}
\vartheta_3(\tau)=\Theta_{\mathbb{Z}}(\tau)=1+2q^1+2q^4+2q^9+\cdots.
\end{equation}
Similarily, the theta series of the $k$-dimensional lattice $\mathbb{Z}^k$ is then
\begin{equation}\label{ZkandkLambda}
\Theta_{\mathbb{Z}^k}(\tau)=\vartheta_3(\tau)^k.
\end{equation}

Theta series of lattices are well studied object in analytic number theory. Here is a well known result concerning theta series of unimodular lattices.
\begin{lemma}(Hecke)\cite{sphere packings} \label{Hecke}
If $\Lambda$ is a unimodular lattice then
$$\Theta_{\Lambda}(\tau)\in \mathbb{C}[\vartheta_3(\tau),\Delta_8(\tau)].$$
\end{lemma}

This lemma tells us that the theta series of any unimodular lattice can be generated by $\vartheta_3(\tau)$ and $\Delta_8(\tau)$. These objects we have discussed in the last part of this subsection are actually \textit{modular forms}. Interested readers may refer to \cite{modular forms} for an introduction.

\subsection{Previous results}

Lattice encoding for the wiretap channel (\ref{channel model}) is done via a generic coset coding strategy \cite{ISITA}: let $\Lambda_e\subset\Lambda_b$ be two nested lattices. A $k$-bit message is mapped to a coset in $\Lambda_b/\Lambda_e$, after which a vector is randomly chosen from the coset as the encoded word. The lattice $\Lambda_e$ can be interpreted as introducing confusion for Eve, while $\Lambda_b$ is intended to ensure reliability for Bob. Since a message is now corresponding to a coset of codewords instead of one single codeword, the probability of correct decoding is then summing over the whole coset (suppose that we do not have power constraint and are utilizing the whole lattice to do the encoding). Here we are interested in computing Eve's probability of correct decision 
$$
P_{c,e}= \sum_{\mathbf{t}\in \Lambda_e}\int_{\mathcal{V}_{\Lambda_b(\mathbf{x}+\mathbf{t})}}\frac{1}{(\sigma_e\sqrt{2\pi})^n}e^{-||\mathbf{y}-\mathbf{x}||^2/2\sigma_e^2}d\mathbf{y}.
$$ 
With a change of variable and applying the \textit{Taylor expansion} to the exponential function, the value of $P_{c,e}$ is approximated in \cite{Gaussian wiretap codes} by
$$
\sum_{\mathbf{t}\in \Lambda_e}e^{-||\mathbf{t}||^2/2\sigma_e^2}\frac{\sqrt{\mbox{det}\Lambda_b}}{(\sigma_e\sqrt{2\pi})^n}\left(1-\frac{\mathcal{U}(\mathcal{V}_{\Lambda_b})}{2\sigma_e^2\sqrt{\mbox{det}\Lambda_b}}\right),
$$
where 
$$
\mathcal{U}(\mathcal{V}_{\Lambda_b})=\int_{\mathcal{V}_{\Lambda_b}}||\mathbf{x}||^2d\mathbf{u}
$$
is the unnormalized second moment of $\Lambda_b$. Since $\sqrt{\mbox{det}\Lambda_b}$ and $\mathcal{U}(\mathcal{V}_{\Lambda_b})$ are invariants of $\Lambda_b$, to minimize $P_{c,e}$ is then to minimize 
\begin{equation}\label{secrecy theta series}
\sum_{\mathbf{t}\in \Lambda_e}e^{-||\mathbf{t}||^2/2\sigma_e^2},
\end{equation}
which is easily recognized as the theta series of $\Lambda_e$ at $\tau=\frac{i}{2\pi\sigma_e^2}$. 

Motivated by the above argument, the confusion brought by the lattice $\Lambda_e$ with respect to no coding (namely, use a scaled version of the lattice $\mathbb{Z}^n$ with the same volume) is measured as follows:
\begin{definition} \cite{ISITA}
Let $\Lambda$ be an $n$-dimensional lattice of volume $v^n$. The secrecy function of $\Lambda$ is given by
$$\Xi_{\Lambda}(\tau)=\frac{\Theta_{v\mathbb{Z}^n}(\tau)}{\Theta_{\Lambda}(\tau)}, \tau\in \mathcal{H}.$$
The \textit{secrecy gain} is then the maximal value of the secrecy function with respect to $\tau$ and is denoted by $\chi_{\Lambda}$.
\end{definition}

\begin{figure}[htp]
\label{fig:secrecy function of E_8}
\centering
\includegraphics[width=80mm, height=60mm]{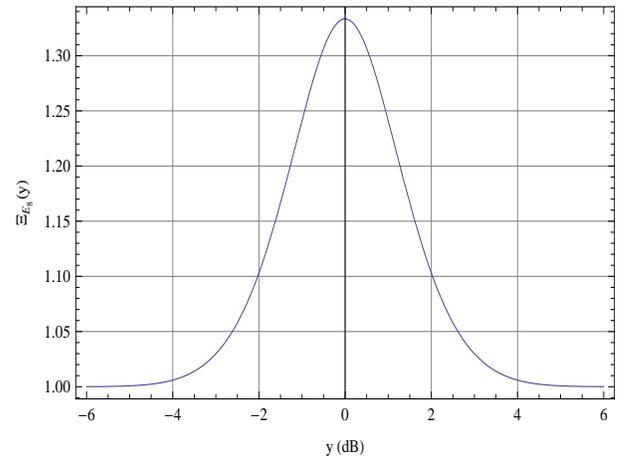}
\caption{Secrecy function of $E_8$}
\end{figure}

A large class of lattices was shown to have a symmetry point (called \textit{weak secrecy gain}) at $\tau=\frac{i}{(\mbox{det}\Lambda)^{\frac{1}{n}}}$ in their secrecy function through the \textit{Poisson summation formula} \cite{Gaussian wiretap codes}. For example, Fig. 1 
shows the secrecy function $\Xi_{E_8}(\tau)$ of $E_8$, where 1) we set $y=-i\tau$ and restrict to real positive values of $y$, since by (\ref{secrecy theta series}) we are only interested in the values of $\Theta_{E_8}(\tau)$ with $\tau=yi,~y>0$ and 2) $y$ is plotted in decibels to transform the multiplicative symmetry point into an additive symmetry point. The symmetry point can be seen to be $y=0$ dB corresponding to $y=1$, and hence to $\tau=i$. This class of lattices contains lattices whose duals are obtained from themselves by possibly a rotation, reflection, and change of scale. Let us now focus on unimodular lattices, for which we have $\Lambda^{*}=\Lambda$ by definition. It was a conjecture by Belfiore and Sol\'e \cite{ITW}, that for these lattices, $\tau=i$ is not only the symmetry point, but also the point achieving the secrecy gain:
\begin{equation}\label{secrecy gain}
\chi_{\Lambda}=\Xi_{\Lambda}(i)=\frac{\vartheta_3^n(i)}{\Theta_{\Lambda}(i)}.
\end{equation}
This conjecture was recently proven by A.-M. Ernvall-Hyt\"onen \cite{EH,EHITW} for a special class of lattices called \textit{extremal even unimodular lattices}. 
The idea of the proof is to write the secrecy function of a lattice $\Lambda$ as a function of the quantity 
$$z=\frac{16\Delta_8(\tau)}{\vartheta_3^8(\tau)}=\frac{\vartheta_2^4(\tau)\vartheta_4^4(\tau)}{\vartheta_3^8(\tau)},\ \tau\in \mathcal{H}.$$ 
She shows that
$$z\in[0,\frac{1}{4}]$$
and that the maximum $\frac{1}{4}$ of $z$ is achieved at $\tau=i$. The rest of the proof consists of showing that the function $f_{\Lambda}(z)$ is increasing in $[0,\frac{1}{4}]$. Later we will prove the conjecture for extremal odd unimodular lattices as well as unimodular lattices in small dimensions using this idea.

\section{The secrecy gain of unimodular lattices}\label{sec:results}
For the sake of convenience, we will assume that the symmetry point $\tau=i$ is really the maximum of the secrecy function through this section. We will then justify the claim for the specific lattices we discuss, but note that the general conjecture is still open.

\subsection{A general formula}

We are now ready to give our first result, namely a general formula for the secrecy gain of unimodular lattices. From Lemma \ref{Hecke} we have the following decomposition of the theta series of a unimodular lattice $\Lambda$:
\begin{equation} \label{unimodular}
\Theta_{\Lambda}(\tau)=\sum_{r=0}^{[\frac{n}{8}]}a_r\vartheta_3^{n-8r}(\tau)\Delta_8^r(\tau), a_r \in \mathbb{Z}.
\end{equation}
Consequently, the reciprocal of the secrecy gain of $\Lambda$ is
$$
\begin{array}{ll}
        1/\chi_{\Lambda}&=\frac{\Theta_\Lambda(i)}{\vartheta_3(i)^n}\\
                                    &=\frac{\sum_{r=0}^{[\frac{n}{8}]}a_r\vartheta_3^{n-8r}(i)\Delta_8^r(i)}{\vartheta_3^n(i)}\\
                                    &=\sum_{r=0}^{[\frac{n}{8}]}a_r(\frac{\Delta_8(i)}{\vartheta_3^8(i)})^r\\
                                   &=\sum_{r=0}^{[\frac{n}{8}]}a_r(\frac{\vartheta_2^4(i)\vartheta_4^4(i)}{16\vartheta_3^8(i)})^r\\
                                   & =\sum_{r=0}^{[\frac{n}{8}]}a_r(\frac{1}{2^6})^r,\\

\end{array}
$$
where the first equality follows from (\ref{secrecy gain}), the second from (\ref{unimodular}), the fourth from (\ref{equ:discriminant}), and the final equality from the following two useful equations concerning the Jacobi theta functions at $\tau=i$ \cite{sphere packings}:
\begin{equation} \label{Jacobi's}
\vartheta_2(i)=\vartheta_4(i) \mbox{ and } \vartheta_3(i)=\sqrt[4]{2}\vartheta_4(i).
\end{equation}

To summarize:
\begin{theorem} The secrecy gain of a unimodular lattice $\Lambda$ of dimension $n$ can be written as
\begin{equation}\label{unimodular secrecy gain}
\chi_{\Lambda}=\frac{1}{\sum_{r=0}^{[\frac{n}{8}]}a_r(\frac{1}{2^6})^r},
\end{equation}
where the $a_i$'s are the coefficients in (\ref{unimodular}).
\end{theorem}

This generalizes the formula for the even case in \cite{Gaussian wiretap codes}.

Fig. 2 gives the plot of the secrecy function of the odd unimodular lattice $D_{12}^+$ mentioned in Example 2.1 (see the paragraph following Fig. 1 for an explanation of the variable $y$). The maximum can be seen to be $\frac{8}{5}$, which we will verify in the next subsection.
\begin{figure}[htp]
\label{fig:secrecy function of D_12}
\centering
\includegraphics[width=80mm, height=60mm]{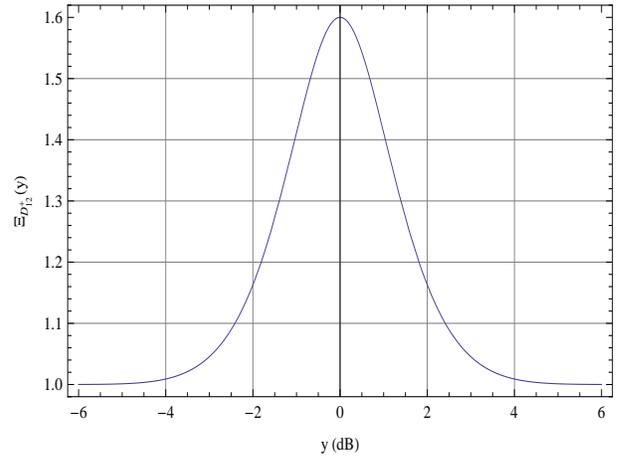}
\caption{Secrecy function of $D_{12}^+$}
\end{figure}

\subsection{Extremal odd unimodular lattices}

In order to find good Gaussian wiretap lattice codes, we look for unimodular lattices with high secrecy gain.
We start by restricting our search to the class of extremal unimodular lattices.

\begin{definition}
Let $\Lambda$ be a lattice of dimension $n$. $\Lambda$ is said to be an extremal lattice if its minimal norm is $[\frac{n}{8}]+1$.\footnote{The definition of extremal has changed. Here we use the earlier version.}
\end{definition}


By definition, an extremal unimodular lattice $\Lambda$ of dimension $n$ contains no vector of norm $1, 2, \cdots, [\frac{n}{8}]$, thus
the coefficients of $q,\ q^2,\ \cdots,\ q^{[\frac{n}{8}]}$ in the theta series given in (\ref{equ:theta series}) are all $0$'s. But by expanding (\ref{unimodular}), we can form another formal sum with coefficients represented as linear combinations of $a_i$'s. Then by comparing the first $[\frac{n}{8}]+1$ terms of the two formal sums, we have a system of $[\frac{n}{8}]$ linear equations in $[\frac{n}{8}]$ unknowns $a_1,\ a_2,\ \cdots,\ a_{[\frac{n}{8}]}$ ($a_0$ is obviously $1$), from which a unique solution can be found. In this way, the secrecy gain of each extremal unimodular lattice can be computed. We illustrate this technique by computing the secrecy gain of $D_{12}^+$ and $O_{23}$.

\textbf{Secrecy gain of $D_{12}^+$}. The theta series of $D_{12}^+$ looks like
$$\Theta_{D_{12}^+}(\tau)=1+0q+A_2q^2+\cdots,\ A_2\neq 0.$$
On the other hand, by (\ref{unimodular}), (\ref{equ:theta3}) and (\ref{equ:discriminant}),
$$
\begin{array}{ll}
\Theta_{D_{12}^+}(\tau)&=\vartheta_3^{12}(\tau)+a_1\vartheta_3^4\Delta_8(\tau)\\
                            &=(1+2q+\cdots)^{12}+a_1(1+2q+\cdots)^4(q+\cdots)\\
                            &=(1+24q+\cdots)+a_1(q+\cdots)\\
                            &=1+(24+a_1)q+\cdots.\\
\end{array}
$$
We now have one linear equation in one unknown $a_1$
$$24+a_1=0,$$
which gives $a_1=-24$,
yielding the secrecy gain
\begin{equation}
\chi_{D_{12}^+}=\frac{1}{1-\frac{24}{2^6}}=\frac{8}{5}.
\end{equation}

\textbf{Secrecy gain of $O_{23}$}. The  theta series of the \textit{Shorter Leech Lattice} $O_{23}$ again looks like
$$\Theta_{O_{23}}(\tau)=1+0q+0q^2+A_3q^3+\cdots,\ A_3\neq 0.$$
On the other hand, by (\ref{unimodular}), (\ref{equ:theta3}) and (\ref{equ:discriminant})
$$
\begin{array}{ll}
\Theta_{O_{23}}(\tau)&=\vartheta_3^{23}(\tau)+a_1\vartheta_3^{15}(\tau)\Delta_8(\tau)+a_2\vartheta_3^7(\tau)\Delta_8^2(\tau)\\
                            &=(1+2q+\cdots)^{23}\\
                            &\ \ +a_1(1+2q+\cdots)^{15}(q-8q^2+\cdots)\\
                            &\ \ +a_2(1+2q+\cdots)^7(q-8q^2+\cdots)^2\\
                            &=(1+46q+1012q^2+\cdots)\\
                            &\ \ +a_1(q+22q^2+\cdots)+a_2(q^2+\cdots)\\
                            &=1+(46+a_1)q\\
                            &\ \ +(1012+22a_1+a_2)q^2+\cdots.\\
\end{array}
$$
This time, we have two linear equations in $a_1$ and $a_2$
$$\left\{ \begin{array}{cc}
46+a_1&=0\\
1012+22a_1+a_2&=0\\
               \end{array}
\right.
,$$
which gives $a_1=-46$ and $a_2=0$,
yielding
\begin{equation}\label{O_23secrecygain}
\chi_{O_{23}}=\frac{1}{1-\frac{46}{2^6}}=\frac{32}{9}.
\end{equation}


By applying this method, we have computed the secrecy gain for each extremal odd unimodular lattice in dimension 
$n,\ n\geq 10$ (see \cite{sphere packings} for a classification), as shown in
Table I. 
A similar table for the even case can be found in \cite{ITW,Gaussian wiretap codes}.

\begin{table}\label{table:extremal odd}
\label{table:secrecy-gains-of-extremal-odd-unimodular-lattices}
\centering
\caption{secrecy gains of extremal odd unimodular lattices}
\begin{tabular}{|c|c|c|c|}
\hline
dim & lattice                          & theta series                                              &secrecy gain \\
\hline
\hline
$12$  & $D_{12}^+$   &$\vartheta_3^{12}-24\vartheta_3^4\Delta_8$  &$\frac{8}{5}$\\
\hline
\hline
$14$  & $(E_7^2)^+$  &$\vartheta_3^{14}-28\vartheta_3^6\Delta_8$ &$\frac{16}{9}$\\
\hline
\hline
$15$  & $A_{15}^+$   &$\vartheta_3^{15}-30\vartheta_3^7\Delta_8 $&$\frac{32}{17}$\\
\hline
\hline
$23$  & $O_{23}$     &$\vartheta_3^{23}-46\vartheta_3^{15}\Delta_8$    &$\frac{32}{9}$\\
\hline
\end{tabular}

\end{table}


\begin{proposition}\label{odd extremal}
The secrecy conjecture is true, namely,  the maximum of the secrecy function is achieved at $\tau=i$ for extremal odd unimodular lattices and the secrecy gains are given as in Table I. 
\end{proposition}

\textit{Proof.} The secrecy gains are computed as illustrated in the examples of $D_{12}^+$ and $O_{23}$. Now we only need to show that the secrecy gains of these unimodular lattices are indeed achieved at $\tau=i$. Recall the definition of secrecy function and the theta series of $D_{12}^+$, which we have just computed. We have that 
$$
\begin{array}{ll}
\Xi_{D_{12}^+}(\tau)&=\frac{\vartheta_3^{12}(\tau)}{\vartheta_3^{12}(\tau)-24\vartheta_3^4(\tau)\Delta_8(\tau)}\\
                                &=\frac{1}{1-\frac{24z}{16}},\\
\end{array}
$$
where $z=\frac{\vartheta_2^4(\tau)\vartheta_4^4(\tau)}{\vartheta_3^8(\tau)}$. Recall also that it was shown in \cite{EH} that $z\in[0,\frac{1}{4}]$ and $\frac{1}{4}$ is achieved at $\tau=i$. It then suffices to show that the denominator is decreasing in $[0,\frac{1}{4}]$, which is obviously true, since its derivative is negative in $[0,\frac{1}{4}]$. Thus the maximum of the secrecy function is achieved at $z=\frac{1}{4}$, namely, $\tau=i$. 

We do the same for the other three extremal odd unimodular lattices, namely:
$$
\begin{array}{ll}
\Xi_{(E_7^2)^+}(\tau)&=\frac{\vartheta_3^{14}(\tau)}{\vartheta_3^{14}(\tau)-28\vartheta_3^6(\tau)\Delta_8(\tau)}\\
                                &=\frac{1}{1-\frac{28z}{16}},\\
\end{array}
$$
$$
\begin{array}{ll}
\Xi_{A_{15}^+}(\tau)&=\frac{\vartheta_3^{15}(\tau)}{\vartheta_3^{15}(\tau)-30\vartheta_3^7(\tau)\Delta_8(\tau)}\\
                                &=\frac{1}{1-\frac{30z}{16}},\\
\end{array}
$$
$$
\begin{array}{ll}
\Xi_{O_{23}}(\tau)&=\frac{\vartheta_3^{23}(\tau)}{\vartheta_3^{23}(\tau)-46\vartheta_3^{15}(\tau)\Delta_8(\tau)}\\
                                &=\frac{1}{1-\frac{46z}{16}}.\\
\end{array}
$$
That the maximum of each secrecy function is achieved at $z=\frac{1}{4}$, namely, $\tau=i$ follows similarly. The proof is completed.
\bigskip

A unimodular lattice containing vectors of norm $1$ can  always be written as the direct sum of a unimodular lattice without vectors of norm $1$ and a cubic lattice $\mathbb{Z}^k$ \cite{sphere packings}. From the definition of the secrecy function, we have that the secrecy gain is determined by the component that contains no vector of norm $1$. In fact,
$$\chi_{\Lambda\oplus \mathbb{Z}^k}=\frac{\vartheta_3^n(i)}{\Theta_{\Lambda}(i)\vartheta_3^k(i)}=\frac{\vartheta_3^{n-k}(i)}{\Theta_{\Lambda}(i)}=\chi_{\Lambda}.$$

By refering to the enumeration of unimodular lattices \cite{sphere packings}, the lattices $E_8$, $D_{12}^+$, $(E_7^2)^+$ and $A_{15}^+$ are the only unimodular lattices that do not contain vectors of norm $1$ in dimensions less than $16$. The secrecy gain of $E_8$ was computed in \cite{ITW,Gaussian wiretap codes} and observe that the secrecy gains of these three lattices are already given in Table I. Thus we in fact have all the best unimodular lattices in dimension $n$, $8<n<16$, namely, $E_8\oplus \mathbb{Z}$ in dimension $9$, $E_8\oplus \mathbb{Z}^2$ in dimension $10$, $E_8\oplus \mathbb{Z}^3$ in dimension $11$, $D_{12}^+$ in dimension $12$, $D_{12}^+\oplus \mathbb{Z}$ in dimension $13$,  $(E_7^2)^+$ in dimension $14$ and $A_{15}^+$ in dimension $15$.

We will deal with unimodular lattices in dimension $n$, $16\leq n\leq 23$ in the next subsection.

\subsection{Unimodular lattices in small dimensions}

The computation of secrecy gain of extremal unimodular lattices can easily be adapted to cover a large family of unimodular lattices, namely, non-extremal unimodular lattices that do not contain vectors of norm $1$ in dimensions $16\leq n \leq 23$. 


We show the computation of the secrecy gain of $(D_8^2)^+$ to illustrate the technique before deriving a general formula of the secrecy gain for all the $111$ lattices and proving the secrecy conjecture for these lattices. 

\textbf{Secrecy gain of $(D_8^2)^+$}. The lattice $(D_8^2)^+$ does not contain any vector of norm $1$. Thus the corresponding coefficient $A_1$ in the theta series is $0$. Its kissing number is $224$, which means that the first nonzero coefficient $A_2=224$ and the theta series of $(D_8^2)^+$ looks like
\begin{equation}\label{Ai's}
\Theta_{(D_8^2)^+}(\tau)=1+0q+224q^2+A_3q^3+\cdots.
\end{equation}
On the other hand, by (\ref{unimodular}), (\ref{equ:theta3}) and (\ref{equ:discriminant}),
$$
\begin{array}{ll}
\Theta_{(D_8^2)^+}(\tau)&=\vartheta_3^{16}(\tau)+a_1\vartheta_3^8(\tau)\Delta_8(\tau)+a_2\Delta_8^2(\tau)\\
                            &=(1+2q+\cdots)^{16}\\
                            &\ \ +a_1(1+2q+\cdots)^8(q-8q^2+\cdots)\\
                            &\ \ +a_2(q-8q^2+\cdots)^2\\
                            &=(1+32q+480q^2+\cdots)\\
                            &\ \ +a_1(q+8q^2+\cdots)+a_2(q^2+\cdots)\\
                            &=1+(32+a_1)q+(480+8a_1+a_2)q^2+\cdots.\\
\end{array}
$$
This time, we have two linear equations in $a_1$ and $a_2$
$$\left\{ \begin{array}{ccc}
32+a_1&=&0\\
480+8a_1+a_2&=&224\\
               \end{array}
\right.
,$$
which gives $a_1=-32$ and $a_2=0$,
yielding
\begin{equation}
\chi_{(D_8^2)^+}=\frac{1}{1-\frac{32}{2^6}}=2.
\end{equation}


We now derive a general formula for the secrecy gain of all the $111$ non-extremal unimodular lattices. 
\begin{proposition}\label{small dimensions}
The secrecy gain conjecture is true, namely, the secrecy gain is achieved at $\tau=i$ for non-extremal unimodular lattices in dimension $16\leq n \leq 23$ and the secrecy gain is given by
\begin{equation}\label{secrecygaingeneral}
\chi_{\Lambda}=\frac{1}{1-\frac{2n}{2^6}+\frac{2n(n-23)+K(\Lambda)}{2^{12}}},
\end{equation}
where $K(\Lambda)$ denotes the kissing number of $\Lambda$.
\end{proposition}

\textit{Proof.} The theta series of a lattice $\Lambda$ in question looks like
\begin{equation}\label{Ai'sformula}
\Theta_{\Lambda}(\tau)=1+0q+K(\Lambda)q^2+A_3q^3+\cdots.
\end{equation}
On the other hand, by (\ref{unimodular}), (\ref{equ:theta3}) and (\ref{equ:discriminant})
$$
\begin{array}{ll}
                       
\Theta_{\Lambda}(\tau)&=\vartheta_3^{n}(\tau)+a_1\vartheta_3^{n-8}(\tau)\Delta_8(\tau)+a_2\vartheta_3^{n-16}(\tau)\Delta_8^2(\tau)\\
                            &=(1+2q+\cdots)^{n}\\
                            &\ \ +a_1(1+2q+\cdots)^{n-8}(q-8q^2+\cdots)\\
                            &\ \ +a_2(1+2q+\cdots)^{n-16}(q-8q^2+\cdots)^2\\
                            &=(1+2\left( \begin{aligned} n\\ 1\\\end{aligned} \right)q+2^2\left( \begin{aligned} n\\ 2\\\end{aligned} \right)q^2+\cdots)\\
                            &\ \ +a_1(1+2\left( \begin{array}{c} n-8\\ 1\\\end{array} \right)q+\cdots)(q-8q^2+\cdots)\\
                            &\ \ +a_2(q^2+\cdots)\\
                            &=(1+2nq+2n(n-1)q^2+\cdots)\\
                            &\ \ +a_1(q+(2n-24)q^2+\cdots)\\
                            &\ \ +a_2(q^2+\cdots)\\
                            &=1+(2n+a_1)q\\       
                            &\ \ +(2n(n-1)+(2n-24)a_1+a_2)q^2+\cdots.\\
\end{array}
$$
Now by comparing the two expressions of $\Theta_{\Lambda}$, we have two linear equations in $a_1$ and $a_2$
$$\left\{ \begin{array}{ccc}
2n+a_1&=&0\\
2n(n-1)+(2n-24)a_1+a_2&=&K(\Lambda)\\
               \end{array}
\right.
,$$
which gives $a_1=-2n$ and $a_2=2n(n-23)+K(\Lambda)$, yielding from the conjecture
$$
\chi_{\Lambda}=\frac{1}{1-\frac{2n}{2^6}+\frac{2n(n-23)+K(\Lambda)}{2^{12}}}.
$$


We have yet to show that the maximum is indeed achieved at $\tau=i$. Recalling the definition of secrecy function and the theta series we have just computed, the secrecy function of $\Lambda$ can be written as
$$
\begin{array}{ll}
\Xi_{\Lambda}(\tau)&=\frac{\vartheta_3^n(\tau)}{\vartheta_3^n(\tau)-2n\vartheta_3^{n-8}(\tau)\Delta_8(\tau)+(2n(n-23)+K(\Lambda))\vartheta_3^{n-16}(\tau)\Delta_8^2(\tau)}\\
                               &=\frac{1}{1-\frac{2n}{16}z+\frac{2n(n-23)+K(\Lambda)}{16^2}z^2}\\
                               &=\frac{1}{D(z)},\\
\end{array}
$$
where $D(z)=1-\frac{2n}{16}z+\frac{2n(n-23)+K(\Lambda)}{16^2}z^2$ and $z=\frac{\vartheta_2^4(\tau)\vartheta_4^4(\tau)}{\vartheta_3^8(\tau)}$. Recall again that it was shown in \cite{EH} that $z\in[0,\frac{1}{4}]$ and $\frac{1}{4}$ is achieved at $\tau=i$.  It suffices to show that the denominator $D(z)$ is decreasing in $[0,\frac{1}{4}]$. We now examine the derivative of the denominator. Note that $z\in[0,\frac{1}{4}]$, $16\leq n\leq 23$ and the largest kissing number for these lattice is $760$. Now,
$$
\begin{array}{ll}
D^{'}(z)&=-\frac{n}{8}+\frac{2n(n-23)+K(\Lambda)}{2^7}z\\
&\leq-\frac{n}{8}+\frac{2n(n-23)+K(\Lambda)}{2^9}\\
&=\frac{-64n+2n(n-23)+K(\Lambda)}{2^9}\\
&<\frac{-1024+0+K(\Lambda)}{2^9}\\
&=\frac{K(\Lambda)-1024}{2^9}\\
&\leq \frac{760-1024}{2^9}\\
&<0.\\
\end{array}
$$
This tells us that the denominator $D(z)$ is decreasing in $[0,\frac{1}{4}]$ and the maximum of the secrecy function is achieved at $z=\frac{1}{4}$, namely, $\tau=i$. The proof is completed.
\bigskip

\begin{table}
\label{table:secrecy gains of 16-23}
\centering
\caption{Secrecy gains of non-extremal unimodular lattices of  dimension $n$, $16\leq n\leq 23$}
\begin{tabular}{|c|c|c|c|}
\hline
dim  &lattice                                  &kissing number      &secrecy gain\\
\hline
\hline
$16$&$E_8^2$                               &$480$                        &$\frac{16}{9}$\\ 
\hline
$16$&$D_{16}$                              &$480$                        &$\frac{16}{9}$\\   
\hline  
\boldmath{$16$}&\boldmath{$D_8^2$}                              &\boldmath{$224$}                        &\boldmath{$2$}\\   
\hline
\hline
\boldmath{$17$}&\boldmath{$A_{11}E_6$}                         &\boldmath{$204$}                        &\boldmath{$\frac{32}{15}$}\\   
\hline
\hline
$18$&$A_{17}A_1$                        &$308$                        &$\frac{32}{15}$\\  
\hline
$18$&$D_{10}E_7A_1$                  &$308$                        &$\frac{32}{15}$\\  
\hline
\boldmath{$18$}&\boldmath{$D_6^3$}                             &\boldmath{$180$}                        &\boldmath{$\frac{16}{7}$}\\  
\hline
\boldmath{$18$}&\boldmath{$A_9^2$}                             &\boldmath{$180$}                        &\boldmath{$\frac{16}{7}$}\\  
\hline
\hline
$19$&$E_6^3O_1$                       &$216$                        &$\frac{64}{27}$\\  
\hline
$19$&$A_{11}D_7O_1$                &$216$                        &$\frac{64}{27}$\\  
\hline
\boldmath{$19$}&\boldmath{$A_7^2D_5$}                      &\boldmath{$152$}                        &\boldmath{$\frac{32}{13}$}\\  
\hline
\hline
$20$&$D_20$                              &$760$                        &$\frac{32}{17}$\\  
\hline
$20$&$D_{12}E_8$                      &$504$                        &$\frac{32}{15}$\\  
\hline
$20$&$D_{12}D_8$                     &$376$                        &$\frac{16}{7}$\\  
\hline
$20$&$E_7^2D_6$                      &$312$                        &$\frac{64}{27}$\\  
\hline
$20$&$A_{15}D_5$                     &$280$                        &$\frac{128}{53}$\\  
\hline
$20$&$D_8^2D_4$                     &$248$                        &$\frac{32}{13}$\\  
\hline
$20$&$A_{11}E_6A_3$                &$216$                        &$\frac{128}{51}$\\  
\hline
$20$&$D_6^3A_1^2$                 &$184$                        &$\frac{64}{25}$\\  
\hline
$20$&$A_9^2A_1^2$                 &$184$                        &$\frac{64}{25}$\\  
\hline
$20$&$A_7^2D_5O_1$               &$152$                        &$\frac{128}{49}$\\  
\hline
\boldmath{$20$}&\boldmath{$D_4^5$}                           &\boldmath{$120$}                        &\boldmath{$\frac{8}{3}$}\\  
\hline
\boldmath{$20$}&\boldmath{$A_5^4$}                           &\boldmath{$120$}                        &\boldmath{$\frac{8}{3}$}\\  
\hline
\hline
$21$&$A_{20}O_1$                     &$420$                       &$\frac{256}{109}$\\  
\hline
$21$&$A_{13}E_7O_1$               &$308$                        &$\frac{128}{51}$\\  
\hline
$21$&$A_{11}D_9O_1$              &$276$                        &$\frac{64}{25}$\\  
\hline
$21$&$A_{12}A_8O_1$              &$228$                         &$\frac{256}{97}$\\  
\hline
$21$&$D_7A_7E_6O_1$            &$212$                         &$\frac{8}{3}$\\  
\hline
$21$&$A_9D_6A_5O_1$            &$180$                        &$\frac{128}{47}$\\  
\hline
$21$&$A_8^3A_4O_1$              &$164$                        &$\frac{256}{93}$\\  
\hline
$21$&$A_7D_5^2A_3O_1$       &$148$                         &$\frac{64}{23}$\\  
\hline
$21$&$A_6^3A_2O_1$             &$132$                          &$\frac{256}{91}$\\  
\hline
$21$&$A_5^3D_4A_1O_1$       &$116$                          &$\frac{128}{45}$\\  
\hline
$21$&$A_4^5O_1$                   &$100$                          &$\frac{256}{89}$\\  
\hline
\boldmath{$21$}&\boldmath{$A_3^7$}                          &\boldmath{$84$}                           &\boldmath{$\frac{32}{11}$}\\  
\hline
\hline
$22$&$D_{14}E_7A_1$              &$492$        &$\frac{64}{27}$\\
\hline
$22$&$E_8E_7^2$                     &$492$        &$\frac{64}{27}$\\
\hline
$22$&$D_{10}^2A_1^2$            &$364$        &$\frac{64}{25}$\\
\hline
$22$&$A_{15}D_6O_1$              &$300$        &$\frac{8}{3}$\\
\hline
$22$&$D_{10}D_6^2$                &$300$        &$\frac{8}{3}$\\
\hline
$22$&$D_8E_7D_6A_1$            &$300$        &$\frac{8}{3}$\\
\hline
$22$&$A_{13}D_7A_1O_1$        &$268$        &$\frac{128}{47}$\\
\hline
$22$&$D_8D_6^2A_1^2$          &$236$        &$\frac{64}{23}$\\
\hline
$22$&$A_{10}^2O_2$                &$220$        &$\frac{256}{91}$\\
\hline
$22$&$E_6^2A_5^2$                 &$204$        &$\frac{128}{45}$\\
\hline
$22$&$A_{11}D_5A_5A_1$        &$204$        &$\frac{128}{45}$\\
\hline
$22$&$A_9D_7A_5O_1$           &$204$        &$\frac{128}{45}$\\
\hline
$22$&$A_9E_6D_5A_1O_1$      &$204$        &$\frac{128}{45}$\\
\hline
$22$&$D_6^2D_4^2A_1^2$     &$172$        &$\frac{32}{11}$\\
\hline
$22$&$A_7^2D_6O_2$             &$172$        &$\frac{32}{11}$\\
\hline
$22$&$A_9A_7D_4A_1O_2$     &$172$        &$\frac{32}{11}$\\
\hline
$22$&$A_8A_6^2O_2$             &$156$        &$\frac{256}{87}$\\
\hline
$22$&$A_7^2A_3^2A_1^2$      &$140$        &$\frac{128}{43}$\\
\hline
$22$&$D_5^2A_5^2O_2$         &$140$        &$\frac{128}{43}$\\
\hline
$22$&$A_7D_5A_5A_3A_1O_1$&$140$        &$\frac{128}{43}$\\
\hline
$22$&$A_6^2A_4^2O_2$         &$124$        &$\frac{256}{85}$\\
\hline
$22$&$D_4^4A_1^6$               &$108$        &$\frac{64}{21}$\\
\hline
$22$&$A_5^3A_3A_1^3O_1$   &$108$        &$\frac{64}{21}$\\
\hline
$22$&$A_5^2D_4A_3^2O_2$   &$108$        &$\frac{64}{21}$\\
\hline
\end{tabular}
\end{table}

\begin{table}
\centering
\begin{tabular}{|c|c|c|c|}
\hline
dim  &lattice                                  &kissing number      &secrecy gain\\
\hline
$22$&$A_4^4A_2^2O_2$         &$92$        &$\frac{256}{83}$\\
\hline
$22$&$A_3^6A_1^2O_2$         &$76$         &$\frac{128}{41}$\\
\hline
$22$&$A_2^{10}O_2$               &$60$        &$\frac{256}{81}$\\
\hline
\boldmath{$22$}&\boldmath{$A_1^{22}$}                     &\boldmath{$44$}        &\boldmath{$\frac{16}{5}$}\\
\hline
\hline
$23$&$A_{15}E_8$                   &$480$        &$\frac{128}{51}$\\
\hline
$23$&$A_{19}A_4$                   &$400$        &$\frac{256}{97}$\\
\hline
$23$&$D_{11}A_{11}O_1$        &$352$        &$\frac{128}{47}$\\
\hline
$23$&$A_{11}E_7A_5$             &$288$        &$\frac{128}{45}$\\
\hline
$23$&$A_9E_7E_6O_1$           &$288$        &$\frac{128}{45}$\\
\hline
$23$&$D_9E_6^2O_2$            &$288$        &$\frac{128}{45}$\\
\hline
$23$&$A_{14}E_6A_2O_1$      &$288$        &$\frac{128}{45}$\\
\hline
$23$&$D_9A_7^2$                  &$256$        &$\frac{32}{11}$\\
\hline
$23$&$A_{13}A_8A_1O_1$      &$256$        &$\frac{32}{11}$\\
\hline
$23$&$A_{11}D_8A_3O_1$      &$256$        &$\frac{32}{11}$\\
\hline
$23$&$D_8A_7^2O_1$            &$224$        &$\frac{128}{43}$\\
\hline
$23$&$D_7^2A_7O_2$            &$224$        &$\frac{128}{43}$\\
\hline
$23$&$A_{11}A_7A_4O_1$      &$208$        &$\frac{256}{85}$\\
\hline
$23$&$A_{10}A_9A_2A_1O_1$&$208$        &$\frac{256}{85}$\\
\hline
$23$&$E_6D_5^3O_2$            &$192$        &$\frac{64}{21}$\\
\hline
$23$&$E_6D_6A_5^2O_1$      &$192$        &$\frac{64}{21}$\\
\hline
$23$&$D_7A_7D_5A_3O_1$   &$192$        &$\frac{64}{21}$\\
\hline
$23$&$A_8E_6A_6A_2O_1$    &$192$        &$\frac{64}{21}$\\
\hline
$23$&$A_{10}A_6D_5O_2$     &$192$        &$\frac{64}{21}$\\
\hline
$23$&$A_9D_6D_5A_1O_2$   &$192$        &$\frac{64}{21}$\\
\hline
$23$&$D_5^4O_3$                 &$160$        &$\frac{128}{41}$\\
\hline
$23$&$A_9A_5A_4^2O_1$     &$160$        &$\frac{128}{41}$\\
\hline
$23$&$D_6D_5A_5^2O_2$     &$160$        &$\frac{128}{41}$\\
\hline
$23$&$A_7D_6A_5A_3A_1O_1$&$160$      &$\frac{128}{41}$\\
\hline
$23$&$A_8A_6D_5A_2O_2$    &$160$        &$\frac{128}{41}$\\
\hline
$23$&$A_8A_7A_5A_1O_2$    &$160$        &$\frac{128}{41}$\\
\hline
$23$&$A_8A_5^2A_2^2O_1$  &$144$        &$\frac{256}{81}$\\
\hline
$23$&$A_7^2A_4A_3A_2$      &$144$        &$\frac{256}{81}$\\
\hline
$23$&$A_7A_6^2A_1^2O_2$  &$144$        &$\frac{256}{81}$\\
\hline
$23$&$D_5^2A_3^4O_1$        &$128$        &$\frac{16}{5}$\\
\hline
$23$&$A_7D_4^2A_3^2O_2$  &$128$        &$\frac{16}{5}$\\
\hline
$23$&$A_7A_5A_4^2A_1O_2$ &$128$        &$\frac{16}{5}$\\
\hline
$23$&$D_5A_5^2D_4A_1^2$   &$128$        &$\frac{16}{5}$\\
\hline
$23$&$A_6^2D_4A_4O_3$       &$128$        &$\frac{16}{5}$\\
\hline
$23$&$A_6D_5A_4^2A_2O_2$ &$128$        &$\frac{16}{5}$\\
\hline
$23$&$A_5^3A_4A_1O_3$        &$112$        &$\frac{256}{79}$\\
\hline
$23$&$A_6A_5A_4A_3A_2A_1O_2$&$112$  &$\frac{256}{79}$\\
\hline
$23$&$D_4^2A_3^4O_3$          &$96$          &$\frac{128}{39}$\\
\hline
$23$&$A_5^2A_3A_2^4O_2$     &$96$         &$\frac{128}{39}$\\
\hline
$23$&$A_5D_4A_3^3A_1^3O_2$&$96$        &$\frac{128}{39}$\\
\hline
$23$&$D_4A_4^3A_2^2O_3$      &$96$        &$\frac{128}{39}$\\
\hline
$23$&$A_5A_4^2A_3^2A_1O_3$&$96$        &$\frac{128}{39}$\\
\hline
$23$&$A_4A_3^5O_4$                &$80$        &$\frac{256}{77}$\\
\hline
$23$&$A_4^2A_3^2A_2^2A_1^2O_3$&$80$&$\frac{256}{77}$\\
\hline
$23$&$A_3^4A_1^8O_3$             &$64$        &$\frac{64}{19}$\\
\hline
$23$&$A_3^3A_2^4A_1^2O_4$   &$64$        &$\frac{64}{19}$\\
\hline
$23$&$A_2^6A_1^6O_5$             &$48$        &$\frac{256}{75}$\\
\hline
\boldmath{$23$}&\boldmath{$A_1^{16}O_7$}                   &\boldmath{$32$}        &\boldmath{$\frac{128}{37}$}\\
\hline
\end{tabular}
\end{table}

Table II \footnote{$O_k$ in the table denotes an empty component of dimension $k$, namely, one containing no vector of norm less than or equal to $2$. Also, for the sake of simplicity, we omit the ``$^+$" which denotes the existence of glue vectors.} 
summarizes the secrecy gains we have computed. Observe that
\begin{enumerate}
\item In dimension $16$, the odd unimodular lattice $(D_8^2)^+$ has secrecy gain $2$, which outperforms its two even counterparts $(E_8^2)^+$ and $D_{16}^+$, both with secrecy gain $\frac{16}{9}$. \\
\item In fact, when the dimension $n$ is fixed the secrecy gain is totally determined by the kissing number $A_2$. The lattice with the best secrecy gain (in boldface) is the one with the smallest kissing number, which can also be seen directly from (\ref{secrecygaingeneral}). This agrees with the observation in \cite{Gaussian wiretap codes} that the best secrecy gain is achieved by extremal lattices, for being extremal in this special case is equivalent to having $A_2=0$. We do not know yet if the secrecy gain is connected to the kissing number in general.
\end{enumerate}

In \cite{Gaussian wiretap codes}, a lower bound on the minimal secrecy gain as a function of $n$ from Siegel-Weil formula for even unimodular lattices was computed. In Fig. 3, the points corresponding to best unimodular lattices are compared to the bound. Note that all the points are the secrecy gains of odd lattices, except for $E_8$ in dimension $8$. We observe that when $n$ grows, the gap between the lower bound and the best lattices decreases, as suggested in \cite{Gaussian wiretap codes}, where it was shown that when $n$ increases, the difference of secrecy gain becomes negligible.
\begin{figure}[htp]
\label{fig:bound}
\centering
\includegraphics[width=80mm, height=60mm]{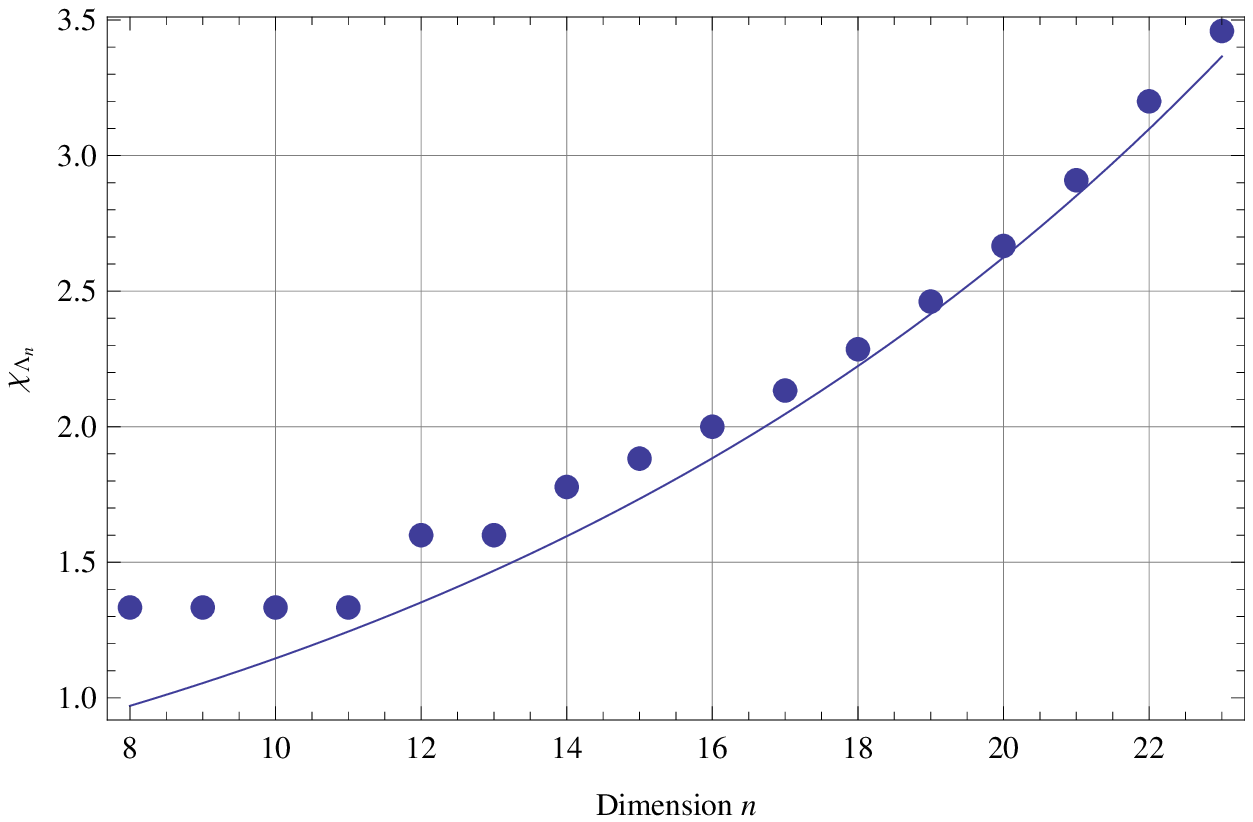}
\caption{Lower bound of the minimal secrecy gain as a function of $n$ from Siegel-Weil formula for even unimodular lattices. Points correspond to best unimodular lattices.}
\end{figure}
\section{Gaussian wiretap codes from unimodular lattices}\label{sec:examples}
As mentioned in Section II, the secrecy gain of a lattice $\Lambda$ characterizes the amount of confusion at Eve that is gained by using this lattice $\Lambda$ as $\Lambda_e$ in the lattice coset code $\Lambda_e\subset \Lambda_b$. Now that we have established the secrecy gain of all the unimodular lattices in dimension smaller than $24$, we need to be able to use these lattices, particularly those with the highest secrecy gain to provide lattice coset codes. To do so, lattice encoding should be performed, which can be handled via Construction A, assuming that we can associate to the chosen lattice a suitable error correction code. We will use some terminology from classical error correction codes in this session. Unfamiliar readers can refer to \cite{error correction codes}.

\subsection{Construction A}

There is a classic way of constructing lattices from binary linear codes called \textit{Construction A}. 
Let $\rho: \mathbb{Z}^n\rightarrow \mathbb{F}_2^n$ be the map of componentwise reduction modulo 2 defined on $\mathbb{Z}^n$.
Let $C$ be a binary $[n,k,d]$ code. Then
$\rho^{-1}(C)$ is a free Abelian group of rank $n$ and hence is a lattice in $\mathbb{R}^n$.

\begin{definition}
The lattice $\Gamma_C$ generated by $C$ is defined by
$$\Gamma_C :=\frac{1}{\sqrt2}\rho^{-1}(C).$$
\end{definition}

To help identify which, if any, error correction code corresponds to a given lattice, we use the following known results:
\begin{theorem}\cite{lattices and codes}\label{unimodular lattices and self-dual codes}
\label{theorem: Construction A}
Let $C$ be a binary linear code and $\Gamma_C$ be the lattice generated by $C$. Then
\begin{enumerate}
  \item $C\subset C^{\perp}$ if and only if $\Gamma_C$ is an integral lattice;
  \item $C$ is doubly even if and only if $\Gamma_C$ is an even lattice;
  \item $C$ is self-dual if and only if $\Gamma_C$ is unimodular.
\end{enumerate}
\end{theorem}

A self-dual code is always an even code. It is called a \textit{type II} code if it is doubly even and \textit{type I} otherwise \cite{error correction codes}. It then follows from Theorem \ref{unimodular lattices and self-dual codes} that $C$ is a type I (respectively type II) code if and only if $\Gamma_C$ is a type I (respectively type II) lattice.

\begin{theorem}\label{kissing number and weight distribution} \cite{sphere packings} Let $C$ be a binary $[n,k,d]$ linear code with weight distribution $W_C(k),\ k=0,1,\cdots,n$. Then the kissing number $K(\Gamma_{C})$ of the lattice $\Gamma_{C}$ generated by $C$ is given by
$$
K(\Gamma_{C})=
\left\{
\begin{array}{rr}
2^dW_{C}(d) &\mbox{if }d<4,\\
2n+16W_{C}(4) &\mbox{if }d=4,\\
2n&\mbox{if }d>4.\\
\end{array}
\right.
$$
\end{theorem}

Theorem \ref{kissing number and weight distribution} gives a way to find the corresponding unimodular lattice for each self-dual code, assuming that there is only one unimodular lattice having the computed kissing number. When we have more than one unimodular lattice with the same kissing number, more considerations are needed to distinguish them. Table III 
gives the list of type I codes of length $n$ ($8<n \leq 23$) \cite{self-dual codes}. According to Theorem \ref{unimodular lattices and self-dual codes}, the lattices generated by these codes are odd unimodular lattices of dimension $n$ ($8<n \leq 23$). The rest of the work consists of finding out the corresponding lattice for each code, through Theorem \ref{kissing number and weight distribution}.

\begin{table}
\label{table:type I codes}
\caption{Type I codes of length $n$ ($8<n \leq 23$)}
\centering
\begin{tabular}{|c|c|}
  \hline
  Codes      &  \begin{tabular}{c|c} \ \ \ \ \ \ \ \ \ \ \  weight \ \ \ distribution\ \ \ \ \ \ \ \ \ \ \ \ \ &  num\\ \end{tabular}\\
  \hline
  \hline
  $[12,6,4]$  &  \begin{tabular}{c|c}  \ \ \ \  \ \ \ \ \ \ \ \   (1,0,\textbf{15},32,15,0,1)\ \ \ \  \ \ \ \  \ \ \ \ \    & 1         \\ \end{tabular}\\
  \hline
  \hline
  $[14,7,4]$ &  \begin{tabular}{c|c}  \ \ \ \  \ \ \ \  \ \ \  (1,0,\textbf{14},49,49,14,0,1)\ \ \ \  \ \ \ \  \ \   & 1         \\ \end{tabular}\\
  \hline
  \hline
  $[16,8,4]$ &  \begin{tabular}{c|c}  \ \ \ \  \ \ \ \  (1,0,\textbf{12},64,102,64,12,0,1)\ \ \ \ \ \ \ \    & 1         \\ \end{tabular}\\
  \hline
  \hline
  $[18,9,4]$  &
                 \begin{tabular}{c|c}
                 (1,0,\textbf{9},75,171,171,75,9,0,1)   &1\\
                 \hline
                \ \ \ \ \  (1,0,17,51,187,187,51,17,0,1)\ \ \ \ \ \    &1\\
                 \end{tabular}\\
  \hline
  \hline

  $[20,10,4]$ &

                \begin{tabular}{c|c}
                (1,0,\textbf{5},80,250,352,250,80,5,0,1)      & 1\\
                \hline
                (1,0,9,72,246,368,246,72,9,0,1)      & 1\\
                \hline
                (1,0,13,64,242,384,242,64,13,0,1)    & 1\\
                \hline
                (1,0,17,56,238,400,238,56,17,0,1)    & 1\\
                \hline
                (1,0,21,48,234,416,234,48,21,0,1)    & 1\\
                \hline
                (1,0,29,32,226,448,226,32,29,0,1)    & 1\\
                \hline
                \ \ \ \ (1,0,45,0,210,512,210,0,45,0,1)\ \ \ \ \        & 1\\
                \end{tabular} \\

  \hline
  \hline
  $[22,11,6]$ &  \begin{tabular}{c|c}   \ (1,0,\textbf{0},77,330,616,616,330,77,0,0,1)  \ \   & 1\\ \end{tabular}\\
  \hline
   $[22,11,4]$ &

                \begin{tabular}{c|c}
                (1,0,4,73,318,628,628,318,73,4,0,1)      & 1\\
                \hline
                (1,0,8,69,306,640,640,306,69,8,0,1)      & 1\\
                \hline
                (1,0,12,65,294,652,652,294,65,12,0,1)    & 1\\
                \hline
                (1,0,16,61,282,664,664,282,61,16,0,1)    & 2\\
                \hline
                (1,0,20,57,270,676,676,270,57,20,0,1)    & 1\\
                \hline
                (1,0,28,49,246,700,700,246,49,28,0,1) & 2\\
                \end{tabular} \\

  \hline
\end{tabular}
\end{table}

\begin{proposition} Best unimodular lattices of  dimension $n$, $8< n\leq 23$ and their corresponding codes are as classified in Table IV. 
\end{proposition}
\textit{Proof.} First by Tables I and II as well as the observations following Propositions \ref{odd extremal} and \ref{small dimensions}, the best unimodular lattices are as shown in Table IV. Now we find their corresponding codes. That the code for the lattice $E_8$ is the $[8,4,4]$ \textit{Extended Hamming code} was mentioned in \cite{Gaussian wiretap codes}. We only need to show the correspondence for all the even dimensions from $12$ to $22$. And since by Theorem \ref{unimodular lattices and self-dual codes}, the lattices generated by the type I codes in Table III 
are odd unimodular lattices, the correspondence can be shown by finding the corresponding lattice for each type I code. For $n=12$, $14$ and $16$, there is only one code of the respective length and only one odd unimodular lattice of the respective dimension, hence it is clear. Let us now deal with the rest of the even dimensions one by one, from $22$ to $18$. Since according to the observation following Proposition \ref{small dimensions}, the best unimodular lattices are those with the smallest kissing numbers, we can directly search for the codes that give the smallest kissing numbers. For $n=22$, $[22,11,6]$ has a minimum distance of $6$, which is greater than $4$. According to Theorem \ref{kissing number and weight distribution}, the generated lattice has kissing number $44$, which is the smallest a type I code of length $22$ can give. By refering to Table II, 
we know that it is the lattice $(A_1^{22})^+$. For $n=20$, $[20,10,4]$ with $W_C(4)=5$ will give the smallest kissing number and applying Theorem \ref{kissing number and weight distribution} again yields $K(\Lambda_C)=120$. There are two odd unimodular lattices both with the same kissing number $120$, namely, $(A_5^4)^+$ and $(D_4^5)^+$. Finally, for $n=18$, we have a similar situation. The code $[18,9,4]$ with $W_C(4)=9$ gives the smallest kissing number $180$ and there are two odd unimodular lattices, $(A_9^2)^+$ and $(D_6^3)^+$, both with the same kissing number $180$.

For the odd dimensions, the unimodular lattices cannot be obtained from Construction A, since the conditions in Theorem \ref{unimodular lattices and self-dual codes} are necessary and sufficient and there does not exist self-dual binary codes of odd length. The proof is completed.
\bigskip

\begin{table}
\label{table:unimodular lattices}
\caption{Best unimodular lattices of  dimension $n$, $8< n\leq 23$ and the corresponding codes}
\centering
\begin{tabular}{|c|c|c|}
\hline
 dim       &lattice                                           &codes \\
\hline
\hline
$9$        &$E_8\oplus\mathbb{Z}$                & $[8,4,4]$\\
\hline
$10$      &$E_8\oplus\mathbb{Z}^2$            & $[8,4,4]$\\
\hline
$11$      &$E_8\oplus\mathbb{Z}^3$            & $[8,4,4]$\\
\hline
$12$      &$D_{12}^+$                                   & $[12,6,4]$\\
\hline
$13$      &$D_{12}^+\oplus\mathbb{Z}$       & $[12,6,4]$\\
\hline
$14$      &$(E_7^2)^+$                                    &$[14,7,4]$ \\
\hline
$15$      &$A_{15}^+$                                   &\\
\hline
$16$      &$(D_8^2)^+$                                   &$[16,8,4]$\\
\hline
$17$      &$(A_{11}E_6)^+$                            &\\
\hline
$18$      &$(A_9^2)^+$ or $(D_6^3)^+$                                  &$[18,9,4]$ with $W_C(4)=9$\\
\hline
$19$      &$(A_7^2D_5)^+$                           &\\
\hline
$20$      &$(A_5^4)^+$ or $(D_4^5)^+$                                  &$[20,10,4]$ with $W_C(4)=5$\\
\hline
$21$      &$(A_3^7)^+$                                    &\\
\hline
$22$      &$(A_1^{22})^+$                                  &$[22,11,6]$ \\
\hline
$23$      & $O_{23}$                                       &\\
\hline

\end{tabular}
\end{table}

\subsection{Coset encoding}
With Construction A, a unimodular lattice $\Lambda$ can be written as 
$$\sqrt{2}\mathbb{Z}^n+\frac{1}{\sqrt{2}}[n, k, d]\mbox{ or }\bigcup_{c_i\in C}\frac{1}{\sqrt{2}}( 2\mathbb{Z}^n+c_i),$$ 
where $C=[n, k, d]$ is the binary linear code that generates $\Lambda$. The encoding is normally done by mapping $k$ bits of information for a codeword of $C$ and $n\lceil\log_2(m)\rceil$ bits of 
information for a bounded set of $\ZZ^n$ around the origin given by $\{ 0,1,\ldots,m-1\}$. In the case of coset encoding 
for a wiretap channel, we can adapt this encoding by setting $\Lambda_e=\Lambda$, and $\Lambda_b=\sqrt{2}\ZZ^n$, in which 
case, $k$ bits of information are indeed used for a codeword of $C$, thus determining a coset, while the other bits are 
either random or least significant. In doing so, we are increasing the confusion at the eavesdropper, however, there is no special coding for Bob.
Let us write
\begin{equation}\label{eq:constrA}
\Lambda_e = \sqrt{2}\ZZ^n + \frac{1}{\sqrt{2}}[n,k,d].
\end{equation}
Since
\[
\ZZ^n = 2\ZZ^n + [n,n,1],
\]
we have that
\[
\frac{1}{\sqrt{2}}\ZZ^n = \sqrt{2}\ZZ^n + \frac{1}{\sqrt{2}}[n,n,1],
\]
which combined with (\ref{eq:constrA}) yields 
\begin{eqnarray*}
\Lambda_e & = & \frac{1}{\sqrt{2}}\ZZ^n  + \frac{1}{\sqrt{2}}[n,n,1]+\frac{1}{\sqrt{2}}[n,k,d]\\
          & = & \frac{1}{\sqrt{2}}\ZZ^n + \frac{1}{\sqrt{2}}C^\dagger 
\end{eqnarray*}
where by definition $[n,k,d]+C^\dagger=[n,n,1]$. Scaling this last equation, we further obtain
\[
2\Lambda_e = \sqrt{2}\ZZ^n + \sqrt{2}C^\dagger 
\]
which together with (\ref{eq:constrA}) gives
\[
\Lambda_e = 2\Lambda_e +\sqrt{2}C^\dagger +  \frac{1}{\sqrt{2}}[n,k,d].
\]
By doing so, we can alternatively choose $\Lambda_b=2\Lambda_e$ instead of $\sqrt{2}\ZZ^n$.

\section{Conclusion and future work}\label{sec:conclusion}

A recent line of work on lattice codes for Gaussian wiretap channels introduced a new lattice invariant called secrecy gain as a code design criterion which captures the confusion that lattice coding can introduce at an eavesdropper. So far, only the secrecy gain of even unimodular lattices was studied.
In this paper, we pursued the study of unimodular lattices by investigating the case of odd unimodular lattices, which exist in greater number and, unlike even lattices, in any dimension. We provided a general formula for the secrecy gain of unimodular lattices in general. We then computed the secrecy gain for odd unimodular lattices, both extremal, and in small dimensions.
As a result, we gave a classification of the best unimodular lattice wiretap codes in small dimensions.

Future work on unimodular wiretap lattice codes concerns the asymptotic behavior of odd unimodular lattices. More generally, it is of interest to generalize the existing work on unimodular lattices to other classes of lattices.

%
%

\section*{Acknowledgment}
The authors would like to thank the reviewers of the ITW version of this paper for their comments.

%
%

%
%






\begin{thebibliography}{1}
%
\bibitem{ISITA}
J.-C. Belfiore and F. Oggier, ``Secrecy gain: a wiretap lattice code design," ISITA 2010.
\url{http://arXiv:1004.4075v2} [cs.IT].
%
\bibitem{Wyner} A. D. Wyner, ``The wire-tap channel," Bell. Syst. Tech. Journal, vol. 54, October 1975.
%
\bibitem{theoretic security}
Y. Liang, H.V. Poor and S. Shamai, ``Information theoretic security," Foundations and Trends in Communications and Information Theory, Vol. 5, Issue 4-5, 2009, Now Publishers.
%


\bibitem{OW-84}
L. H. Ozarow and A. D. Wyner,``Wire-tap channel II,'' {\em Bell Syst. Tech.
Journal}, vol. 63, no. 10, pp. 2135-2157, Dec. 1984.
%
\bibitem{TDCMM-07}
A. Thangaraj, S. Dihidar, A. R. Calderbank, S.W. McLaughlin, and J.-M. Merolla,``Applications of LDPC Codes to the Wiretap Channel,'' {\em IEEE Transactions on Information Theory}, vol. 53, No. 8, Aug. 2007.
%
\bibitem{IEEE Gaussian wiretap channel}
S. K. Leung-Yan-Cheong and M. E. Hellman, ``The Gaussian wire-tap channel", \textit{IEEE} Trans. Inform. Theory, vol. IT-24, no. 4, pp. 451-456, July 1978.
%
\bibitem{LDPC Gaussian wiretap codes}
D. Klinc, J. Ha, S. McLaughlin, J. Barros, and B. Kwak, ``LDPC codes for the Gaussian wiretap channel," in \textit{Proc. ITW}, Oct. 2009.
%
\bibitem{nested codes}
R. Liu, H.V. Poor, P. Spasojevic,  and Y. Liang,  ``Nested codes for secure transmission'',  in Proc. PIMRC, 2008, pp.1-5. 

\bibitem{ITW}
J.-C. Belfiore and P. Sol\'e, ``Unimodular lattices for the Gaussian Wiretap Channel,'' ITW 2010, Dublin. \url{http://arXiv:1007.0449v1} [cs.IT].
%
\bibitem{Gaussian wiretap codes}
F. Oggier, J.-C. Belfiore, and P. Sol\'e, ``Lattice Coding for the Wiretap Gaussian Channel", \url{http://arXiv:1103.4086v1} [cs.IT], 21 Mar 2011.
%
\bibitem{lattices and codes}
W. Ebeling, ``Lattices and Codes", Advanced Lectures in Mathematics, Vieweg \& Sohn, Verlagsgesellschaft mbH, Braunschweig/Wiesbaden, 1994.

%
\bibitem{sphere packings}
J.H. Conway, N.J.A. Sloane, ``Sphere packings, Lattices and Groups", Third edition, Springer-Verlag, New York, 1998.
%
\bibitem{modular forms}
N. Koblitz, ``Introduction to Elliptic Curves and Modular Forms", Graduate Texts in Math. No. 97, Springer-Verlag, New York, Second edition, 1993.
%
%
\bibitem{EH}
A.-M. Ernvall-Hyt\"onen, ``On a Conjecture by Belfiore and Sol\'e on some Lattices", 
\url{http://arxiv.org/PS\_cache/arxiv/pdf/1104/1104.3739v2.pdf}
%
\bibitem{EHITW} A.-M. Ernvall-Hyt\"onen and C. Hollanti, ``On the Eavesdropper’s Correct Decision in Gaussian and Fading Wiretap Channels Using Lattice Codes," ITW 2011, Paraty. pp. 210-214. 
%

\bibitem{error correction codes} F. J. MacWilliams and N. J. A. Sloane, ``The Theory of Error-Correcting Codes", Amsterdam, The Netherlands: North-Holland, 1977.

%

%
\bibitem{self-dual codes} \url{http://www.cs.umanitoba.ca/~umbilou1/SelfDualCodes/toc.html}.
%
%
\bibitem{ITWversion} F. Lin and F. Oggier, ``Secrecy gain of Gaussian wiretap codes from unimodular lattices," ITW 2011, Paraty. pp. 718-722. 
%
\end{thebibliography}
\end{document}